\documentclass[12pt]{amsart}
\usepackage{amsmath,amsthm,amsfonts,amssymb,latexsym,enumerate,xcolor,mathtools} 
\usepackage{float}
\usepackage[pagebackref]{hyperref} 

\newcommand{\F}{\mathbb{F}}
\newcommand{\N}{\mathbb{N}}
\newcommand{\Q}{\mathbb{Q}}

\newcommand{\C}{\mathbb{C}}

\newcommand{\Syl}{\operatorname{Syl}}

\newcommand{\Aut}{\operatorname{Aut}}

\newcommand{\Irr}{\operatorname{Irr}}
\newcommand{\IBr}{\operatorname{IBr}}

\newcommand{\SL}{\operatorname{SL}}

\newcommand{\Gal}{\operatorname{Gal}}

\newcommand{\Dim}{\operatorname{Dim}}

\newtheorem{thm}{Theorem}[section]
\newtheorem{lem}[thm]{Lemma}

\newtheorem*{thmA}{Theorem A}

\newtheorem*{thmB}{Theorem B}

 \newtheorem*{conA}{Conjecture A}
  \newtheorem*{conB}{Conjecture B}

\theoremstyle{definition}

\newtheorem{exmp}[thm]{Example}
\newtheorem{note}[thm]{Note}
\newtheorem{claim}[thm]{Claim}

\numberwithin{equation}{section}

\def\lc{\left\lfloor}   
\def\rc{\right\rfloor}

\begin{document}


\title[Groups with at most $2$ characters with the same degree]{The Baby Monster is the largest group with at most $2$ irreducible characters with the same degree}

\author{Juan Mart\'inez Madrid}
\address{Departament de Matem\`atiques, Universitat de Val\`encia, 46100
  Burjassot, Val\`encia, Spain}
\email{Juan.Martinez-Madrid@uv.es}

\thanks{Research  supported by Ministerio de Ciencia e Innovaci\'on (Grant PID2022-137612NB-I00 funded by MCIN/AEI/10.13039/501100011033 and ``ERDF A way of making Europe"),   and by Generalitat Valenciana CIAICO/2021/163 and  CIACIF/2021/228. }

\keywords{Brauer's Problem 1, Character degrees}

\subjclass[2020]{Primary 20C15}

\begin{abstract}
We classify all finite groups such that all irreducible character degrees appear with multiplicity at most $2$. As a consequence, we prove that the largest group with at most $2$ irreducible characters of the same degree is the Baby Monster.
\end{abstract}

\maketitle

\pagestyle{myheadings}


\section{Introduction}\label{Section1}

Brauer's problem 1 \cite{Brauer} asks for the classification of the complex group algebras. Let $G$ be a finite group and let $\Irr(G)=\{\chi_1,\ldots ,\chi_k\}$. We define the degree pattern of $G$ as the multiset $\{\chi_1(1),\ldots ,\chi_k(1)\}$.  We know that $\C G\cong \bigoplus_{i=1}^{k}M_{\chi_i(1)}(\C)$.  Hence,  Brauer's problem 1 is equivalent to classifying all possible degree patterns of finite groups.

Even though the complete classification of all degree patterns seems impossible to reach, there are some restrictions on the possible degree patterns of finite groups. For example, if $\{n_1,\ldots ,n_k\}$ is the degree pattern of a finite group $G$, then the number of $1$'s in this pattern is $|G:G'|$ which divides $|G|=\sum_{i=1}^{k}n_i$. In addition, we also have that each $n_i$ is the degree of an irreducible character of $G$ and hence it must divide $\sum_{i=1}^{k}n_i$. We are interested in less trivial restrictions.

Landau's Theorem \cite{landau} asserts that the size of a group can be bounded in terms of the number of conjugacy classes. Thus, there exist finitely many degree patterns of a given length. Conjecture A of \cite{Orig} asked for a general version of  Landau's Theorem. Let $G$ be a group. We define 
$$m(G)=\max_{d\in \N}|\{\chi \in \Irr(G)|\chi(1)=d\}|,$$
that is the maximum number of characters with the same degree.  Conjecture A of \cite{Orig} asked  whether there exists $b : \N \rightarrow \N$ such that $|G|\leq b(m(G))$ for every finite group $G$. By Theorem A of \cite{Orig}  this conjecture is true if and only if it is true  for  symmetric groups. Finally, Theorem 1.2 of \cite{Craven} proved that this conjecture is true for symmetric groups.

Those papers proved that $m(G)$ tends to infinity as $|G|$ grows, but did not provide explicit bounds for $|G|$ in terms of $m(G)$. It could be possible to have explicit values by following the proofs in \cite{Orig} and \cite{Craven}. We include a pair of families of groups in which the bounds are computed.

In the nilpotent case, it is possible to obtain an explicit good bound for $|G|$ in terms of $m(G)$. A.  Jaikin-Zapirain \cite{Andrei} proved that if $G$ is a nilpotent group of size at least $8$, then $G$ contains at least $$\frac{C\log(\log(|G|))\log(|G|)}{\log(\log(\log(|G|)))}$$ conjugacy classes, where $C$ is an explicitly computable constant. Since $\chi(1)$ divides $|G|$ for every $\chi \in \Irr(G)$, and $|G|$ has at most $\log(|G|)$ different divisors, then applying the pigeonhole principle, we deduce that 
$$m(G)>\frac{C\log(\log(|G|))}{\log(\log(\log(|G|)))},$$
which tends to infinity as $|G|$ grows. 

In the case of symmetric groups,  D. Craven  \cite{Craven2}  proved that there exists a constant $a$ such that  
$$\log(m(S_n))\leq a\sqrt{\log(n!)/\log(\log(n!))}-\log(\log(n!))$$
for sufficiently large $n$.

Once known that $|G|$ is $m(G)$-bounded, it is natural to ask which are the groups with $m(G)=1,2$ or $3$. It is not completely trivial to show that the unique group with $m(G)=1$ is the trivial group (see Theorem \ref{casem1} below).

In private communication, A. Moret\'o has conjectured the following.

\begin{conA}
Let $B$ be the Baby Monster group. If $G$ is a group  with $m(G)=2$, then $|G|\leq|B|$.
\end{conA}

\begin{conB}
Let $M$ be the Monster group. If $G$ is a group  with $m(G)=3$, then $|G|\leq|M|$.
\end{conB}

Our main result answers affirmatively Conjecture A.

\begin{thmA}
Let $G$ be a group with $m(G)=2$. Then $|G|\leq|B|$.
\end{thmA}

In fact, Theorem A will  follow from the complete classification of all groups with $m(G)=2$.

\begin{thmB}
Let $G$ be a finite group with $m(G)=2$. Then $G$ is one of the following: 
$$C_2,\Sigma_3, D_{10}, \Sigma_4, A_5, \Sigma_5, A_6,A_7,\Sigma_8, A_9,A_{10}, \Sigma_{10}, A_{16}, L_2(11), L_2(7), L_3(3).2, U_3(5).2$$ 
$$ M_{12}.2, M_{22}, McL, Th, J_2, J_2.2, F_3+, Co_1,B$$
$$ \SL(2,5), (C_2)^4\rtimes A_6, (C_2)^4\rtimes A_7, (C_{11})^2\rtimes \SL(2,5).$$
\end{thmB}

We  divide the proof of Theorem B in four steps.

\begin{itemize}

\item  \underline{First step:} Classify all solvable groups with $m(G)=2$.

\item  \underline{Second step:} Classify all almost-simple groups with $m(G)=2$.

\item  \underline{Third step:} Classify all groups with a unique non-abelian composition factor such that $m(G)=2$.

\item  \underline{Fourth step:} Prove that all groups with $m(G)=2$ have at most one non-abelian composition factor.

\end{itemize}

One could ask if the classification of groups with $m(G)=2$ could be deduced by applying the    techniques  in  \cite{Orig} and \cite{Craven}. Of course, we will apply many ideas from those papers, however, the bounds obtained there are too rough, even for $m(G)=2$.

\section{Structure of groups with $m(G)=2$}

In this section we study the possible composition factors and prime divisors of groups  with $m(G)=2$. We begin by relating $m(G/N)$ with $m(G)$ for a normal subgroup $N$ of $G$.

\begin{lem}\label{quotients}
Let $G$ be a finite group and let $N \trianglelefteq G$. Then $m(G/N)\leq m(G)$.
\begin{proof}
The result follows from the fact that $\Irr(G/N)\subseteq \Irr(G)$.
\end{proof}
\end{lem}

We recall some concepts about fields of values of characters. Let $\chi$ be a character (not necessarily irreducible) of $G$. We define the field of values of $\chi$ as
$$\Q(\chi)=\Q(\chi(g)|g \in G)$$
and we will write $\Q(\chi)$ to denote it. We have that $\Q(\chi)\subseteq \Q_{|G|}$, where $\Q_{|G|}$ denotes the cyclotomic extension of order $|G|$.  Given $\sigma \in \Gal(\Q(\chi)/\Q)$ we define the character $\chi^{\sigma}$ as $\chi^{\sigma}(g)=(\chi(g))^{\sigma}$. We observe that $\chi^{\sigma}$ is irreducible provided that $\chi$ is irreducible.  Now, we relate the fields of values of the irreducible characters with $m(G)$.

\begin{lem}\label{fields}
Let $G$ be a finite group and let $\chi \in \Irr(G)$. Then $|\Q(\chi):\Q|\leq m(G)$.
\begin{proof}
We observe that if $\sigma \in \Gal(\Q(\chi)/\Q)$, then $\chi^{\sigma}\in \Irr(G)$ and $\chi^{\sigma}(1)=\chi(1)$. Thus, $|\Q(\chi):\Q|=|\Gal(\Q(\chi)/\Q)|\leq m(G)$.
\end{proof}
\end{lem}

We recall that a group is said to be rational if $\Q(\chi)=\Q$ for every $\chi\in \Irr(G)$. Analogously, we say that a group is quadratic rational if $|\Q(\chi):\Q|\leq 2$  for every $\chi\in \Irr(G)$. With this, we can classify all groups with $m(G)=1$.

 \begin{thm}\label{casem1}
 Let $G$ be a group with $m(G)=1$. Then $G=1$.
 \begin{proof}
 Assume that $G$ is a non-trivial group with $m(G)=1$. Since $|G:G'|$ is the number of irreducible characters of degree $1$, we deduce that $|G:G'|=1$. Thus, $G$ is a perfect group and hence there exists $N \trianglelefteq G$ such that $G/N$ is non-abelian simple. By Lemma \ref{quotients}, we have that  $1\leq m(G/N)\leq m(G)\leq 1$ and hence $m(G/N)=1$.
 
By Lemma \ref{fields}, we have that $|\Q(\chi):\Q|=1$ for all $\chi \in \Irr(G/N)$. Therefore, $G/N$ is a simple group with  $\Q(\chi)=\Q$ for all  $\chi\in \Irr(G/N)$. Thus, applying Corollary B.1 of  \cite{FeitSeitz}, we have that $G/N=Sp_6(2)$ or $G=O_{8}^{+}(2)'$. However, $m(Sp_6(2))=3$ and $m(O_{8}^{+}(2)')=4$. In both cases, we have a contradiction. 

Thus $G=1$.  
 \end{proof}
 \end{thm}
 
We observe that this result can also be deduced from the main result of \cite{BCH}. This result classifies the groups in which the non-linear characters have pairwise different degrees. In particular, all these groups are solvable. If $G$ is a group with $m(G)=1$, then the degrees of the  non-linear characters are pairwise different. Thus $G$ is solvable, but since $m(G)=1$, we have that $|G:G'|=1$, which implies $G=1$.

 In the case when $m(G)=2$, we have that $G$ is   quadratic rational. In particular, we can use the following result for solvable groups.

\begin{thm}[Theorem A of \cite{Tent}]\label{Tent}
Let $G$ be a quadratic rational solvable group. If $p$ divides $|G|$, then $p \in\{2,3,5,7,13\}$.
\end{thm}

Moreover, in the non-solvable case, the possible non-abelian and non-alternating composition factors of quadratic rational groups are classified by the following result. 

\begin{thm}[Theorem 1.2 of   \cite{Trefethen}]\label{Trefethen}
Let $G$ be a quadratic rational group and let $S$ be a non-abelian composition factor of $G$. Then one of the following holds.
\begin{itemize}
\item [(i)] $S=A_n$ for some $n \geq 5$.

\item [(ii)] $S$ is one of the groups in the list $F$, where $F$ is an explicit and finite list of sporadic and groups of Lie type.
\end{itemize}

\end{thm}

This, together with the following two results will provide a classification of the possible non-abelian composition factors of groups with $m(G)=2$.

\begin{lem}\label{cf}
Let $G$ be a finite group with $m(G)=n$ and let $S$ be a non-abelian composition factor of $G$. Then $S$ cannot possess $n+1$ irreducible characters of the same degree which are extendible to $\Aut(S)$.
\begin{proof}
Let $N/M$ be a chief factor isomorphic to $S^k$ for some $k\geq 1$. Taking quotients, we may assume that $M=1$ and that $C_G(N)=1$. Thus, we may  identity $G$ with a subgroup of $\Gamma=\Aut(S)\wr \Sigma_k$ for some $k$. Assume that there exists $\chi_1,\chi_2,\ldots, \chi_{n+1}\in \Irr(S)$ such that $\chi_1(1)=\ldots=\chi_{n+1}(1)$ and each $\chi_i$ is extendible to $\tilde{\chi}_i\in \Irr(\Aut(S))$.

Let $\Phi_i=\chi_i\times \ldots \times \chi_i\in \Irr(S^k)$ for $i \in \{1,\ldots,n+1\}$. Then $\Phi_i$ is extendible to $\tilde{\Phi}_i=\tilde{\chi}_i\times \ldots \times \tilde{\chi}_i\in \Irr(\Aut(S)^k)$. Moreover, the characters $\tilde{\Phi}_i$ are $\Gamma$-invariant and hence they are extendibles to $\Gamma$, by Lemma 1.3 of \cite{Mattarei}. Therefore, $\Phi_1,\ldots,\Phi_{n+1}$ are $n+1$ characters of the same degree which are extendibles to $\Gamma$ and hence to $G$. Thus, $m(G)\geq n+1$ which is impossible.
\end{proof}
\end{lem}

As in Conjecture A of \cite{Orig}, the case of symmetric groups of our problem is challenging. The following result is due to D. Craven \cite{Craven3}.

\begin{thm}[Craven]\label{Auxiliar}
Let $n\geq 13$ and $n \not=16$. Then there exist $3$ characters in $A_n$ with the same degree which are extendible to $\Sigma_n$.
\end{thm}

It is worth to mention that in \cite{Craven2} Craven had claimed, without proof, that $m(S_n)\geq 3$ for $n\geq 200$. We thank him for doing the remaining cases $n<200$ and for writing the details of the proof of the whole result for $n \geq 200$.


\begin{thm}\label{main}
Let $G$ be a group with $m(G)=2$. If $S$ is a non-abelian composition factor, then either $S=A_{t}$,  for $t \in \{5,\ldots, 12\}\cup \{16\}$ or $S$ is one of the following groups:

$$ M_{12}, M_{22}, J_{2}, Co_{1}, Th$$ $$Fi_{24}', HS, McL, He, Suz, HN, B$$ $$L_2(7), L_2(8), L_2(11), L_2(16), L_2(27), L_3(3), L_3(4), L_4(3),U_3(3), $$ $$U_3(4), U_3(5), U_3(8), U_4(2), U_4(3), U_5(2), U_5(4), U_6(2), S_4(4),  S_6(3),  $$ $$O_7(3), O_{8}^{+}(2),  O_{8}^{+}(3),  O_{8}^{-}(2),  {}^{2}_{}E_{6}(2), F_4(2),{}^{2}_{}F_{4}(2)' , G_2(3), G_2(4),{}^{3}_{}D_{4}(2).$$

\begin{proof}
Assume first that $S=A_t$ for $t\geq 7$. In this case, $\Aut(S)=\Sigma_t$, and the result follows from Theorems \ref{cf} and \ref{Auxiliar}.

Assume now that $S$ is a non-alternating composition factor of $G$. Since $m(G)=2$, we have that $G$ is a quadratic rational group by Lemma \ref{fields}. Thus, $S$ must be one of the groups listed in the list F of  Theorem \ref{Trefethen}. By Lemma \ref{cf}, we have that $S$ cannot possess $3$ characters of the same degree which are extendible to $\Aut(S)$. However, by looking at the ATLAS (see \cite{ATLAS}), we have that  the groups listed in the statement are the only groups in Trefethen's list which satisfy this condition.
\end{proof}
\end{thm}

\section{The solvable case}

In this section, we prove the solvable case of our classification. If $G$ is a solvable group with $m(G)=2$, then $1<|G:G'|=\{\chi\in \Irr(G)|\chi(1)=1 \}\leq m(G)=2$ and hence $|G:G'|=2$. In particular, the unique abelian group with $m(G)=2$ is $C_2$.

\begin{thm}\label{casoab}
Let $G$ be a non-abelian  group with $m(G)=2$ and $G''=1$. Then $G\in \{\Sigma_3,D_{10}\}$.
\begin{proof}
We have that $G/G'=C_2$. Let $\lambda \in \Irr(G')\setminus\{1_{G'}\}$. Since $G'$ is abelian, we have that $\lambda(1)=1$ and, in particular, $\lambda$ cannot be extended to $G$. Since $G/G'$ is cyclic, applying Theorem 11.22 of \cite{Isaacscar}, we have that $\lambda$ can be extended to an irreducible character of $I_{G}(\lambda)$. Since $\lambda$ cannot be extended to $G$,  we have that $\lambda$ cannot be invariant and hence $G'\leq I_{G}(\lambda)<G$. Thus, $G'= I_{G}(\lambda)$ and hence, by the Clifford's correspondence  (see Theorem 6.11 of \cite{Isaacscar}) $\lambda^{G}\in \Irr(G)$. Therefore, $\chi(1)=2$ for every  $\chi \in \Irr(G|G')$. Since $m(G)=2$, we deduce that $|G|=2+a2^2$ for $a\in \{1,2\}$. The unique groups of order $6$ and $10$ with $m(G)=2$ are $\Sigma_3$ and $D_{10}$.
\end{proof}
\end{thm}

\begin{thm}\label{secondcase}
Let $G$ be a  solvable group with $m(G)=2$, $G''>1$ and $G'''=1$. Then $G=\Sigma_4$.
\begin{proof}
Clearly, $G/G''$ is a group with $m(G/G'')=2$ and $(G/G'')''=1$. By Theorem \ref{casoab}, $G/G''\in \{\Sigma_3,D_{10}\}$. By hypothesis, $G''$ is a  non-trivial abelian group.

Let $\lambda \in \Irr(G'')\setminus\{1_{G''}\}$. Since $G/G''\in \{\Sigma_3,D_{10}\}$, we deduce that $G'/G''$ is cyclic of prime order and hence, reasoning as before, we have that $\psi=\lambda^{G'}\in \Irr(G')$. Now, we have two possible options. The first one is that $\psi$ extends to $\chi \in \Irr(G)$ and the second one is that $\chi=\psi^{G}\in \Irr(G)$. In the first case $\chi(1)=|G':G''|$ and in the second case $\chi(1)=|G:G''|$.  Thus, if $\chi \in \Irr(G|G'')$ then $\chi(1)\in \{|G:G''|,|G':G''|\}$. Thus, $|G|=|G:G''|+a|G:G''|^2+b|G':G''|^2$, where $a,b \in \{0,1,2\}$ and $(a,b)\not=(0,0)$.

Now, if $G/G''=\Sigma_3$, then $|G|\in \{15,24,32,51,60,78,81,90\}$ and if $G/G''=D_{10}$ then $|G|\in \{35,60,110,135,160,210,260\}$. Now, the calculations in GAP \cite{gap} show that the only group with $m(G)=2$ whose order lies in any of these sets is $\Sigma_4$.
\end{proof}
\end{thm}

Now, to complete the classification of solvable groups with $m(G)=2$ it suffices to prove that $|G'''|=1$ for all solvable groups with $m(G)=2$. In order to prove it, we need the following result, which we will also need in the non-solvable case.

\begin{lem}\label{MinNorm}
Let $G$ be a group and let $V\trianglelefteq G$ be a non-central solvable minimal normal subgroup of $G$. If $|V|=p^d$ for a prime $p$, then $d=k\cdot n$, where $G/V$ possesses  $k$ irreducible $p$-Brauer characters of degree $n$. 
\begin{proof}
 We know that $V$ is a $p$-elementary abelian group for a prime $p$. Thus,  $V$ can be viewed as an irreducible $\F_{p}[G/V]$-module of dimension $d$, where $\F_{p}$ denotes the field of $p$ elements.  Let $\mathfrak{X}$ be the $\F_p$-representation of $G/V$ given by $V$. We can extend $\mathfrak{X}$ to a representation $\mathfrak{X}^{\overline{\F}_p}$, where $\overline{\F}_p$ is the algebraic closure of $\F_p$.  Since $\overline{\F}_p$ is algebraically closed, it is an splitting field for $G/V$ by Corollary 9.4 of \cite{Isaacscar}. Thus, by Theorem 9.21 of \cite{Isaacscar} each irreducible constituent of  $\mathfrak{X}^{\overline{\F}_p}$ occurs with multiplicity $1$.

Let $\chi$ be the $p$-Brauer  character afforded by $\mathfrak{X}^{\overline{\F}_p}$. Let $\phi, \psi\in \IBr(G/V)$ be irreducible constituents of $\chi$. By Lemma 3 of  \cite{Tent}, there exists $\sigma\in \Gal(\Q_{|G/V|}/\Q)$  such that  $\psi=\phi^{\sigma}$. In particular, $\psi(1)=\phi(1)^{\sigma}=\phi(1)$. Moreover, by the previous comments, we have that each constituent of $\chi$ occurs with multiplicity $1$. 
\end{proof}
\end{lem}

Now, we can finish the classification of solvable groups with $m(G)=2$.

\begin{thm}\label{SolvableCase}
Let $G$ be a solvable group with $m(G)=2$. Then $G\in\{C_2,\Sigma_3,D_{10}, \Sigma_4\}$.
\begin{proof}
 We observe that, by Theorems \ref{casoab} and \ref{secondcase} the result  will follow once we prove  that $G'''=1$. Suppose that $G'''>1$. By Theorem \ref{secondcase}, we have that $G/G'''=\Sigma_4$.  Taking an appropriate quotient, we may assume that $G'''$ is a minimal normal subgroup and hence it is $p$-elementary abelian for some prime $p$.

We claim that $p\not=2$. Assume, by contrary, that  $p=2$.  If $G'''$ is central in $G$, then $|G'''|=2$ and hence $|G|=48$. However, there is no group of order $48$ with $m(G)=2$. Thus, we may assume that $G'''$ is not central in $G$. Since $\IBr(\Sigma_4)$ contains a unique non-trivial $2$-Brauer character of degree $2$, we deduce that $|G'''|=2^2$, by Lemma \ref{MinNorm}. However, here is no group of order $96$ with $m(G)=2$. Thus, the claim follows.

Since $G'''$ is normal and abelian, then $\chi(1)$ divides $|G:G'''|=24$ for all $\chi \in \Irr(G)$ by Theorem 6.15 of \cite{Isaacscar}. Thus, we have that $|G|=24+a2^2+b4^2+c6^2+d8^2+e12^2+f24^2$ where $a\in \{0,1\}$ and  $b,c,d,e,f\in \{0,1,2\}$. By Theorem \ref{Tent} and the previous claim, we also know  that $|G|/24$ must be a power of a prime lying in the set $\{3,5,7,13\}$. Therefore, the only possibilities are $|G| \in  \{120, 168,  312, 600,1176 \}$. However, there is no group with $m(G)=2$ whose order lies in this set.

It follows that $G'''=1$ and the result follows.
\end{proof}
\end{thm}

\section{Groups with a unique non-abelian composition factor}\label{Sec4}

In this section, we classify all groups with $m(G)=2$ and a unique non-abelian composition factor. Let us sketch the proof of this classification.

 We begin by an almost simple group $H$ with $m(H)=2$. Then, we will classify all groups $G$ such that $m(G)=2$ and  $G$ possesses a normal solvable subgroup $N$ such that $G/N\cong H$. To do so, we will begin by finding all groups $L$ with a solvable minimal normal subgroup $V$ such that $L/V\cong H$ and $m(L)=2$. Then, we will repeat the same process with $L$ instead of $H$, and we  will iterate this process until we do not find new groups.

 Note that if  $V$ is a  solvable minimal normal subgroup of $L$, then $|V|=p^f$ for a prime $p$ and one of the following holds
 
 \begin{itemize}
\item [(i)] $V$ is central in $L$.

\item[(ii)]  $V$ is an irreducible $\F_{p}[L/V]$-module.
\end{itemize}

\subsection{Almost simple case}

We begin by classifying  all almost-simple groups with $m(G)=2$.  Let $G$ be an almost simple group with $S \leq G \leq \Aut(S)$. If $m(G)=2$, then $S$ must be one of the groups in Theorem \ref{main}. Thus, we only have to find all groups $G$ with $S \leq G \leq \Aut(S)$ for $S$ one of the groups in \ref{main}. Using GAP, we deduce  that the almost-simple groups with $m(G)=2$ are the following.

\begin{equation} \label{lista1}
  \begin{multlined}
A_5, \Sigma_5, A_6,A_7,\Sigma_8, A_9,A_{10}, \Sigma_{10}, A_{16}, L_2(11), L_2(7), L_3(3).2,\\ U_3(5).2,
M_{12}.2, M_{22}, McL, Th, J_2, J_2.2, F_3+, Co1,B
\end{multlined}
\end{equation}

\subsection{Almost simple groups over minimal normal subgroups}

Now, we classify all groups $G$ such that $m(G)=2$, and $G$ possesses a solvable minimal normal subgroup $V$ such that $G/V$ is almost simple. Since $m(G/V)\leq m(G)=2$, we have that  $G/V$ is one of the groups in the list \ref{lista1}.

\begin{itemize}
\item [(i)] $V$ is central in $G$.

\item[(ii)]  $V$ is an irreducible $\F_{p}[G/V]$-module
\end{itemize}

Case (i), we have that $V$ is a minimal normal subgroup and it is central and hence $|V|=p$. By Corollary 5.4 of \cite{Isaacs}, we have that the Sylow $p$-subgroups of $(G/V)'$ cannot be cyclic. In particular, $p$ must divide $|(G/V)'|$.  Thus, we   only have to check the groups $p.G/V$ for $G/V$ any of the groups in the list \ref{lista1} and $p$ dividing $|(G/V)'|$. We have that $m(G)>2$ unless $G=2.A_5=\SL(2,5)$.

The case (ii) can be studied in two different ways. The first one can be applied for all groups, but it provides rough bounds for $|V|$. On the other hand, the second method provides much better bounds, but it requires an additional hypothesis.

\begin{note}[First way to proceed in case (ii)]\label{FirstMethod}
By Lemma \ref{MinNorm}, we have that $|V|\geq p^{t}$, where $t$ denotes the smallest degree of a non-principal irreducible  $p$-Brauer character.

Since $V$ is normal and abelian, we have that $\chi(1)$ divides $|G/V|$ for every $\chi \in \Irr(G)$. Since $m(G)=2$ we deduce that $|\Irr(G)|\leq 2\cdot d(|G/V|)$ and that for each $d$, dividing $|G|$, there exists at most $2$ irreducible characters of degree $d$. Let $\chi \in \Irr(G)$, let $d=\chi(1)$ and let $\theta \in \Irr(V)$ such that $[\theta, \chi_V]\not=0$. Then, by Clifford's Theorem (see Theorem 6.2 of \cite{Isaacscar})  we have that $|G:I_G(\theta)|$ divides $d$ and there are $|G:I_G(\theta)|$ characters of $V$ lying under $\chi$. Thus, for each character $\chi \in \Irr(G)$ of degree $d$, we obtain at most $d$ irreducible characters in $V$ lying under $\chi$. Therefore,  we deduce that $|V|=|\Irr(V)|\leq 2\cdot (\sum_{d \mid |G/V|}d)$.


From this discussion we deduce that $p^{t}\leq 2\cdot (\sum_{d \mid |G/V|}d) $ and therefore we have that $ t\leq \lc \log_p(2\cdot (\sum_{d \mid |G/V|}d))\rc$.

\end{note}

We show how to apply this method  in an example.

\begin{exmp}[First method in case $G/V\cong B$ and $p=2$.]
Assume that $G/V\cong B$ and that $p=2$. Then the smallest degree of a $2$-Brauer character of $B$ is $4370$ (see Table 1 of \cite{Jansen}). Thus, by the note above, $2^{4370}\leq p^{d}=|V|\leq 2\cdot d(|B|)\cdot |B|$, which is false. Thus, $V=1$.
\end{exmp}

Reasoning similarly, we have that  if the   smallest degree of a non-principal irreducible $p$-Brauer character is larger than $ 2\cdot \lc \log_p(2\cdot (\sum_{d \mid |G/V|}d))\rc$, then $V=1$. This allows us to discard many of the cases. However, in many of our cases, the group $G/V$ has an irreducible $p$-Brauer character with small degree.

For this reason, we introduce a  second method that can only be applied when every irreducible character of $V$ is extendible to its inertia group. This hypothesis will hold in many of the cases that the first method could not discard.

\begin{note}[Second way to proceed in case (ii)]
Let  $G$ be a group with $m(G)=2$ and $V\trianglelefteq G$ be  a $p$-elementary abelian group. Assume that any  $\lambda \in \Irr(V)$ can be extended to an irreducible character in $I_G(\lambda)$.

Let $\lambda \in \Irr(V)$. By hypothesis, we have that $\lambda$ has an extension $\psi\in \Irr(I_G(\lambda))$. Thus, by the Clifford correspondence and Gallagher’s theorem (see Theorem 6.11 and Corollary 6.17 of \cite{Isaacscar}), we have that the  characters $(\psi \theta)^G$ for $\theta \in \Irr(I_G(\lambda)/V)$ are irreducible and different for different $\theta$. It follows that $m(I_G(\lambda)/V)\leq 2$. Assume now that $\lambda_1,\lambda_2\in \Irr(V)$ are not $G$-conjugated and that $I_G(\lambda_1)/V\cong I_G(\lambda_2)/V$. If $m(I_G(\lambda_1)/V)=2$, then we obtain $4$ irreducible characters of the same degree and hence $m(G)>2$.

Thus, we have that
$$|V|\leq 2|G/V|+\sum_{\{H| V<H<G \text{ and }  m(H/V)=2\}}|G:H|+1.$$
More precisely,
$$|V|= a|G/V|+\sum_{\{H| V<H<G \text{ and }  m(H/V)=2\}}a_H|G:H|+1,$$
where $a\in \{0,1,2\}$ and $a_H\in \{0,1\}$.

Observe that the hypotheses of this note hold in the following cases:

\begin{itemize}

\item If $G$ splits over  $V$, then any  $\lambda \in \Irr(V)$ can be extended to a character in $I_G(\lambda)$, by  Exercise 6.18 of \cite{Isaacscar}. Let $H^2(G/V,V)$ be the second cohomology group of $G/V$ on $V$. We recall that the extensions of $G/V$ by $V$ are in one to one correspondence to the cocycles in $H^2(G/V,V)$. As a consequence, if $H^2(G/V,V)=0$, then $G$ splits over $V$.

\item If $P/V\in \Syl_p(G/V) $ and $P/V$ is cyclic (in particular when $|G/V|_p\leq p$), then any  $\lambda \in \Irr(V)$ can be extended to a character in $I_G(\lambda)$ by combining Corollary 11.22 and Theorem 6.26  of \cite{Isaacscar}
\end{itemize}
\end{note}

Let us show how to apply this method in an example.

\begin{exmp}[Second method in case $G/V\cong A_5$ and $p\not=2$] 
We have that $|G/V|_p\leq p$.  We know that the unique proper subgroups of $G/V$ with $m(H)=2$ are $\{C_2,\Sigma_3,D_{10}\}$. Thus, $|V|$ can be expressed as 
$$|V|=1+a\cdot 60+b\cdot 30+c\cdot 10+d\cdot 6,$$
where $0\leq a \leq 2$ and $0\leq b,c,d\leq 1$. In particular, this forces $|V|\leq 167$. We observe that any power of $3$ or $5$ can be expressed as above, which forces $p\geq 7$.

 If $p\geq7$, then the $p$-Brauer characters of $A_5$ are exactly the irreducible complex characters of $A_5$. Since the smallest degree of a non-principal irreducible character of $A_5$ is $3$, we deduce that  $  |V|\geq p^3 \geq 7^3=343$. This is impossible since $|V|\geq 167$.

Thus $V=1$.
\end{exmp}

In the case when $G/V$ is $A_n$ or $S_n$ we have that $G/V$ has  irreducible modules with small dimension. Moreover, it is hard to calculate the second cohomology group of these representations for large $n$. That moves us to study the cohomology of  modules of $A_n$ and $S_n$ over $\F_p$ with small degree.

Let $H\in \{A_n,S_n\}$.  Let $M=\F_p \epsilon_1+\ldots  +\F_p \epsilon_{n}$ be the natural permutation $\F_pH$-module,  $I=\langle \epsilon_1-\epsilon_2, \ldots, \epsilon_{n-1}-\epsilon_{n} \rangle$,  $l =\F_p(\epsilon_1+\ldots +\epsilon_{n})$ and $L=I/l$. If $p\nmid n$, then $I$ is an irreducible $H$-module of dimension $n-1$ and if $p\mid n$, then $L$ is an irreducible $H$-module of dimension $n-2$. The main results of \cite{KP} classify the second cohomology groups of these modules. The original version of that paper had a mistake and  needed a \textit{corrigenda}. The   following theorems correspond to the corrected versions of the results in  \cite{KP}.
 
 \begin{thm}[Corollaries 1 and 2 of \cite{KP}]
 Let $n \geq 5$, let $p$ be a prime and let $H\in \{A_n,\Sigma_n\}$. Then $\Dim_{\F_p}(H^2(H,L))\leq 1$ and  one of the following holds
 
 \begin{itemize}
 \item [(i)] $H^2(H,L)=0$.

 \item [(ii)] $p=2$, $n$ is even and $H=\Sigma_n$.

\item [(iii)] $p=3$, $n \in \{5,8,9\}$ and $H=A_n$.

\item[(iv)] $p=5$, $n=5$ and $H=A_n$.
\end{itemize} 
 \end{thm}

 \begin{thm}[Corollaries 5 and 6 of \cite{KP}]
 Let  $n\geq 5$, let $H\in \{A_n,\Sigma_n\}$ and let $p$ be a prime not dividing $n$. Then    one of the following holds
 
 \begin{itemize}
 \item [(i)] $H^2(H,I)=0$.

\item [(ii)] $p=3$, $n\in\{5,8\}$, $H=A_n$ and $\Dim_{\F_3}(H^2(A_n,I))=1$.

\end{itemize} 
 \end{thm}

Let us follow the following steps:

\begin{itemize}

\item [a)] Use the first method to prove that $V=1$ when the smallest degree of a non-principal $p$-Brauer character is larger $\lc \log_p(2\cdot (\sum_{d \mid |G/V|}d))\rc$.  Note that, in the case of sporadic groups, the smallest degrees of $p$-Brauer characters can be found in \cite{Jansen}. 

\item [b)] Use the second method in the cases remaining from a) to prove that $V=1$ when $|G/V|_p\leq p$ or when $H^2(G/V,V)=0$.

\item [c)] In the cases remaining from b), determine the possibilities for $|V|$  by using any of the following methods:
\begin{itemize}

\item Calculate the irreducible $\F_p[G/V]$-representations  using GAP and MAGMA \cite{MAGMA}.

\item Look at the degrees of the irreducible $\F_p[G/V]$-representations in ATLAS online \cite{ATLASonline}.

\item Apply Lemma \ref{MinNorm} together with the degrees of the $p$-Brauer characters.

\end{itemize}

\end{itemize}

Following these steps prove that $V=1$ unless when $(G/V, Dim_{\F_p}(V))$ appears Table \ref{Dimensions2}.

\begin{table}[h]
\caption{Possible  dimensions of $V$}
\begin{center}
\begin{tabular}{ |c|c|c|c|} 
\hline
$G/V$    & 2  &  3 & 5  \\ \hline
$A_6$ & $4$ & $6$& - \\ \hline
$A_7$ & $4$ & -& - \\ \hline
$\Sigma_{8}$ &  $6,8,14$ & -& - \\ \hline
 $A_9$ & - & $7$& - \\ \hline
 $A_{10}$ & $16$ & - & - \\ \hline
 $\Sigma_{10}$ & $8$,$16$ & -& - \\ \hline
  $L_2(7)$ & $3$ & - & - \\ \hline
  $L_3(3).2$ & $12$ & $7$ & - \\ \hline
  $U_3(5).2$ & $20$ & - & $8$ \\ \hline
  $M_{12}.2$ & $10$& - & - \\ \hline
  $M_{22}$ & $10$& - & - \\ \hline
  $J_2.2$ & $12$ & - & - \\ \hline
  $McL$ & $22$ & - & - \\ \hline
  $Co_1$ & $24$ & - & - \\ \hline
\end{tabular}
\label{Dimensions2}
\end{center}
\end{table}

In the remaining cases of Table \ref{Dimensions2}, we have that it is possible to decide whether $m(G)=2$ or $m(G)>2$ by one of the following methods:

\begin{itemize}

\item The group $G$ is perfect and has order smaller than $2 \cdot 10^6$ and hence they are computable by the {\tt PerfectGroup} command in GAP.

\item The character table of $G$ is recorded in the list AllCharacterTableNames() in GAP. This happens for example in the cases $(C_2^{10})\rtimes M_{22}$ and $(C_2^{24})\rtimes Co_{1}$.

\item Construct all extensions of $V$ by $G/V$ in MAGMA  and compute its  character degrees. To do this

\end{itemize}

After applying these methods, we prove  that if $G$ possesses a solvable minimal normal subgroup $V$ such that $G/V$ is almost simple and $V$ is an irreducible $G/V$-module, then $G\in \{(C_2)^4\rtimes A_6,(C_2)^4\rtimes A_7\}$. We observe that these groups are recorded  as  {\tt PerfectGroup}(5760,1) and {\tt PerfectGroup}(40320,2) in the {\tt PerfectGroup} library of GAP. 

\subsection{The remaining cases}

Now, we have to repeat the same argument as before for $G/V\in \{SL(2,5),(C_{2})^4\rtimes A_6, (C_{2})^4\rtimes A_7 \}$, and $V$ a minimal normal subgroup of $G$. As before, there are  two possible options

\begin{itemize}
\item [(i)] $V$ is central in $G$.

\item[(ii)]  $V$ is an irreducible $\F_{p}[G/V]$-module.

\end{itemize}

Assume first that we are in case (i). As before, we have that $G=p.G/V$ for $p$ a prime dividing $|G/V|$. Since $G/V$ is perfect, we have that $G$ must be perfect. Thus, we only have to check the perfect groups of size $p\cdot n$ for $n\in \{120,5760,40320\}$ and $p\mid n$. In any case, we have that $G$ is a perfect group of order smaller than $10^6$ and hence, $G$ can be found in the {\tt PerfectGroup} library  of  GAP. We have that $m(G)>2$ for all these groups.

Assume now  that $G/V\in \{A_6,A_7\}$ and that $V$ is an irreducible $\F_{p}[G/V]$-module. Assume first that $p\geq 5$. In this case, $|G/V|_p\leq p$ and hence we may apply the second method to deduce that $m(G)>2$. Thus, we may assume that $p\in \{2,3\}$.  In all cases, we have that there exists $V\leq U \trianglelefteq G$ such that $G/U\in \{A_6,A_7\}$, $U/V=(C_2)^4$ and $V$ is elementary abelian. Computing the irreducible representations of $G/V$ over $\F_p$ for $p\in \{2,3\}$, we have that if $\Psi$ is the representation of $G/V$ over $\F_p$ afforded by $V$, then $U$ is contained in the kernel of $\Psi$.  If $p=3$ then, we have that $U=W\times V$  and hence $m(G/W)=2$, which implies $V=1$ by the study in the  previous subsection.  Thus, it only remains to study the case when $p=2$.

In any of these cases, $U$ is contained in the Kernel of the representation afforded by $V$  and hence $V$ is an irreducible $\F_2$-module of $G/U$. Moreover, $H^2(G/V,V)=0$ unless when $G/V= (C_{2})^4\rtimes A_7$ and $|V|\in\{2^6,2^{14}\}$ and hence the extension of $G$ splits over $V$. Thus, if $G$ is none of these cases, we may apply the second method to prove that $m(G)>2$. Thus, we may assume that   $G/V= (C_{2})^4\rtimes A_7$ and $|V|\in\{2^6,2^{14}\}$. In these cases, we may construct all possible extensions to prove that $m(G)>2$.

Assume now that  $G/V=\SL(2,5)$ and that $V$  is an irreducible $\F_{p}[\SL(2,5)]$-module.  Assume first that $p=2$. We know that the dimension of the irreducible representations of $\SL(2,5)$ over $\F_2$ are $1$ and $4$. Thus, $|G|=120\cdot 2^4=1920$. Since there is no group of order $1920$ with $m(G)=2$, we deduce that $V=1$.  Assume now that $p>2$. In this case, applying the second method, we  deduce that $|V|=1+a\cdot 120+b\cdot 60$, where $a \in \{0,1,2\}$ and $b \in \{0,1\}$. This forces $|V|=121$ and $G=(C_{11}\times C_{11})\rtimes SL(2,5)$ and we observe that $m(G)=2$. Thus, we will also have to deal with the case $G/V=(C_{11}\times C_{11})\rtimes SL(2,5)$.

\subsection{The last case}

Let us assume that $G/V=(C_{11})^2 \rtimes \SL(2,5)$ and $V$ is a minimal normal subgroup. As before, $V$ is either central in $G$ or it is an irreducible $G/V$-module.

If $V$ is central, then looking at  the perfect groups of order $p\cdot 14520 $  for $p\in \{2,3,5,11\}$, we deduce that $m(G)>2$.

Assume that $V$ is an irreducible $\F_{p}[G/V]$-module. Then, there exists $V\leq U \trianglelefteq G$ such that $G/U= \SL(2,5)$ and $U/V=(C_{11})^2$.   Now, we have that if $\Psi$ is a representation of $G/V$ over $\F_p$, then either $U$ is contained in the Kernel of $\Psi$ or  $p=2$ and $\Psi(1)=120$. Since $2^{120}> 2\cdot (\sum_{d \mid |\SL(2,5)|}d) $, we will assume that $U$ is contained in the Kernel of the representation afforded by $V$. If $p\not=11$ then, we have that $U=W\times V$ and hence $m(G/W)=2$, which implies $V=1$ by  the study in the  previous subsection. Thus, it only remains to prove the case $p=11$. If $|V|\geq 11^3$, then $H^2(G/V,V)=0$  and hence we may apply the second method to prove that $m(G)>2$. Thus, we may assume that $|G|=11^2\cdot 14520=1756920$. By examination of the groups of order $1756920$ in the {\tt PerfectGroup} library of GAP, we have that $m(G)>2$.

\vspace{1cm}

Therefore, summarizing all the results in this section, we have proved the following result.

\begin{thm}\label{1cf}
Let $G$ be a group with $m(G)=2$ and  a unique non-abelian composition factor. Then $G$ is one of the following groups:
$$A_5, \Sigma_5, A_6,A_7,\Sigma_8, A_9,A_{10}, \Sigma_{10}, A_{16}, L_2(11), L_2(7), L_3(3).2, U_3(5).2$$ 
$$ M_{12}.2, M_{22}, McL, Th, J_2, J_2.2, F_3+, Co_1,B$$
$$ \SL(2,5), (C_2)^4\rtimes A_6, (C_2)^4\rtimes A_7, (C_{11})^2\rtimes \SL(2,5).$$
\end{thm}

\section{Proof of Theorem B}

We aim to prove that if $m(G)=2$, then $G$ possesses at most one  non-abelian composition factor. Assume by contrary that there exists a group $G$ with $m(G)=2$ and at least $2$ non-abelian composition factors. We reason towards a contradiction to prove that such a group $G$ does not exists.

 Without loss of generality, we may assume that $G$ is a counterexample of minimal order. Let $N$ be a proper normal subgroup of $G$. Then $m(G/N)\leq m(G)=2$ and $|G/N|<|G|$ and hence, by the minimality of $G$, we deduce that $G/N$ possesses at most one non-abelian composition factor and hence $G/N$ is one of the groups in Theorem \ref{SolvableCase} or in Theorem \ref{1cf}. Moreover, if $N \trianglelefteq G$ is solvable, then $G$ possesses at most one non-abelian composition factor, which is impossible. Thus, the solvable radical of $G$ is trivial and all minimal normal subgroups of $G$ are product of non-abelian simple groups.

 We divide the proof in three different claims.

\begin{claim}[First claim: The group $G$ possesses a unique minimal normal subgroup]
Assume  that there exists $N,M\trianglelefteq G$ two minimal normal subgroups. By the previous comments, we have that both $G/N$ and $G/M$ are groups with at most $1$ non-abelian composition factor. That implies that $N$ and $M$ are simple groups and both $G/N$ and $G/M$ are almost simple groups. Thus, $G/N$ and $G/M$ are one of the almost-simple groups obtained in list \ref{lista1}, $Soc(G/M)=N$ and $Soc(G/N)=M$ and $$N\times M\leq G\leq G/M\times G/N\leq \Aut(N)\times \Aut(M).$$ 
By examination of the groups in list \ref{lista1}, we have that $|G/N|\in \{|M|,2|M|\}$ and hence $|G|\in \{|M||N|,2|M||N|\}$. If $|G|=|M||N|$, then $G=N\times M $ and hence $m(G)=m(N\times M )\geq m(N)m(M)\geq 4$, a contradiction. Thus, $|G|=2|M||N|$. Since $|G:M\times N|=2$ we have that  $m(M\times N)\leq 5$ by using Clifford theory.

Assume  that $G/M$ is simple. In this case, $G/M=N$ and hence $G\leq N \times G/N\leq N\times \Aut(M)$. Let $\chi_1,\chi_2\in \Irr(N)$ with $\chi_1(1)=\chi_2(1)$. Then the characters $\chi_1\times 1_M,\chi_2\times 1_M\in \Irr(N\times M)$ are extendible to  $N\times \Aut(M)$ and hence, they are extendibles to $G$. Let $\Phi_1,\Phi_2\in \Irr(G)$ such that $\Phi_1$ extends $\chi_1\times 1_M$ and $\Phi_2$ extends $\chi_2\times 1_M$.  Since $G/(M\times N)\not=1$, we have that there exist $\lambda, \rho\in \Irr(G/N)$ such that $\lambda(1)=\rho(1)$. Thus, using Gallagher's Theorem, we have that $\lambda\Phi_1,\lambda\Phi_2,\rho\Phi_1,\rho\Phi_2\in \Irr(G)$ are all different and have the same degree, a contradiction.

Thus $G/M$ and $G/N$ are not  simple groups and hence $G/M,G/N\in \{\Sigma_5,\Sigma_8, \Sigma_{10}, \\  J_{2}.2,L_3(3).2, M_{12}.2, U_3(5).2\}$. Assume that $G/M\not \in \{\Sigma_5, \Sigma_{10}, J_{2}.2\}$. In this case, $m(N)=3$ and hence $m(N\times M)\geq  m(M)m(N) \geq 6$. However, since  $|G:M\times N|=2$  and $m(G)>2$, then applying Clifford theory, we have that  $m(M\times N)\leq 5$. Therefore, we deduce that  $G/M, G/N \in \{\Sigma_5, \Sigma_{10}, J_{2}.2\}$. 

Assume that $M=N$. Let $\chi\in \Irr(N)\setminus \{1\}$ be extendible to $\Irr(\Aut(N))$. Then $\chi\times 1_M,1_N\times \chi\in \Irr(M\times N)$ are extendibles to $ \Aut(N)\times \Aut(M)$ and hence are extendibles to $G$. Using Gallagher's Theorem again, we obtain $4$ characters in $G$ with the same degree.

Therefore, we only have to consider the cases $(G/M,G/N)\in \{(\Sigma_{5},\Sigma_{10}),(\Sigma_{5},J_{2}.2),\\(\Sigma_{10},J_{2}.2)\}$. We study each case separately:
\begin{itemize}

\item \underline{Case $(G/M,G/N)=(\Sigma_{10},J_{2}.2)$:} In this case, we have that there exist $\chi \in \Irr(A_{10})$ and $\psi\in \Irr(J_2)$ such that $\chi(1)=\psi(1)=36$ and are extendibles to $\Sigma_{10}$ and $J_{2}.2$ respectively. Thus, $\chi\times 1 ,1\times \psi\in \Irr(A_{10}\times J_2)$ and they are extendibles to  $\Sigma_{10}\times J_{2}.2$ and hence to $G$. Thus, we we obtain $4$ characters in $G$ with the same degree, which is impossible.

\item \underline{Case $(G/M,G/N)=(\Sigma_{5},\Sigma_{10})$:} In this case, we have that there exist $\chi \in \Irr(A_5)$ and $\psi,\eta\in \Irr(A_{10})$ such that $\chi(1)=4$, $\psi(1)=75$, $\eta(1)=300$ and are extendibles to $\Sigma_5$ and $\Sigma_{10}$ respectively. Thus, $\chi\times \psi ,1\times \eta\in \Irr(A_5\times A_{10})$ have the same degree and they are  extendibles to $\Sigma_{5}\times \Sigma_{10}$ and hence to $G$, which is impossible.

\item \underline{Case $(G/M,G/N)=(\Sigma_{5},J_{2}.2)$:} Let $\chi\in \Irr(J_2)$ be extendible to $J_2.2$ with degree $126$. We have that $1 \times \chi \in \Irr(A_{5}\times J_2)$ is extendible to $G$. Thus, by Gallagher's Theorem, $G$ possesses $2$ characters of degree $126$ lying over $1 \times \chi$. We also have that there exist $\chi_1,\chi_2\in \Irr(A_5)$ and  $\psi_1,\psi_2\in \Irr(J_2)$ with $\chi_1(1)=\chi_2(1)=3$  and $\psi_1(1)=\psi_2(1)=21$. Then all characters in the set $\{\chi_i\times\psi_j|1\leq i,j\leq 2\}$ are irreducible in $A_{5}\times J_2$ and have degree $63$. If all characters in this set are extendible to $G$, then $G$ possesses $4$ characters with degree $63$. Thus, one of these characters is not extendible. There is no loss to assume that $\chi_1\times\psi_1$ is not extendible to $G$. Since $G/(N\times M)\cong C_2$, we have that $N\times M$ is the inertia group of $\chi_1\times\psi_1$. Thus, by Clifford's correspondence, we have that $(\chi_1\times\psi_1)^G\in \Irr(G)$ and has degree and does not lie over $1 \times \chi$. Thus, $G$ possesses $3$ characters of degree $126$, which is impossible.
\end{itemize}

Thus, the claim follows.
\end{claim}

Thus, there exists a  unique minimal normal subgroup $N$ in $G$. Since $N$ is non-solvable, we have that $C_G(N)=1$. Thus, if $N$ is the copy of $n$ copies of a non-abelian simple group $S$, then we can embed $G$ into $\Gamma:=\Aut(S^n)\cong\Aut(S)\wr \Sigma_n$. Moreover, if $K=G\cap \Aut(S)^n$, then $K\trianglelefteq G$, $K/N$ is solvable and $G/K$ can be embedded into $\Sigma_n$.

\begin{claim}[Second claim:  $K/N\not=1$]
Assume that $K=N$. Then $G$ can be embedded into $\Delta=S\wr \Sigma_n$. Let $\phi,\psi\in \Irr(S)$ such that $\phi(1)=\psi(1)$ and $\phi\not=\psi$. Let $\Phi=\phi\times \ldots \times \phi$ and $\Psi=\psi\times \ldots \times \psi$ with $n$ copies of $\phi$ and $\psi$, respectively. Thus, we have that $\Phi$ and $\Psi$ are two $\Delta$-invariant characters in of $N$. Thus, by Lemma 1.3 of \cite{Mattarei}, we have that $\Phi$ and $\Psi$ are extendibles to $\Delta$ and hence, they are extendible to $G$. Reasoning as before, we have that $m(G)\geq 4$, a contradiction.
\end{claim}

As a consequence, we have that $G/N$ cannot be an almost-simple group and that $G/N\not=C_2$. Thus, we may assume that $G/N\in \{\Sigma_3,D_{10}, \Sigma_4, \SL(2,5), C_{11}^{2}\rtimes \SL(2,5), C_2^4\rtimes A_6, C_2^4\rtimes A_7\}$ and that $K/N\not=1$.

\begin{claim}[Third claim:  The group $G$ cannot exist]
By the minimality of $G$, we know that $G/N$ is a group with $m(G/N)=2$ and a unique non-abelian composition factor. Thus, $G$ is one of the groups in Theorem \ref{SolvableCase} or in Theorem \ref{1cf}. By the second claim, we have  that $K/N\not=1$.   Since $K/N$ is a normal solvable subgroup of $G/N$, we deduce that $G/N$ cannot be almost simple. As a consequence, we have that $G/N\in \{ \Sigma_3,D_{10}, \Sigma_4, \SL(2,5), C_{11}^{2}\rtimes \SL(2,5), C_2^4\rtimes A_6, C_2^4\rtimes A_7\}$, $G/K\leq \Sigma_n$ and $K/N$ contains the unique normal subgroup in $G/N$.

Let $k=\lc n/2\rc$, let $\phi\in \Irr(S) \setminus \{1_S\}$ be extendible to $\chi \in \Irr(\Aut(S))$ and let $\Phi\in \Irr(S^n)$ be a product of $k$ copies of $\phi$ and $n-k$ copies of $1_S$.  We know that $\Phi$ can be extended to $\Psi\in \Irr(\Aut(S)^n)$ which is the product of  $k$ copies of $\chi$ and $n-k$ copies of $1_S$. By Lemma 1.3 of \cite{Mattarei}, we have that $\Psi$ is extendible to $I_{\Gamma}(\Psi)=I_{\Gamma}(\Phi)$, which is isomorphic to $\Aut(S)^n\rtimes  (\Sigma_k\times \Sigma_{n-k})$. Thus, $\Phi $ is extendible to an irreducible character of its inertia group and since $G\leq \Gamma$, then $\Phi$ is extendible to $I_{\Gamma}(\Phi)\cap G=I_{G}(\Phi)$. By Gallagher's Theorem and Clifford's correspondence, we have that $m(I_{G}(\Phi)/N)\leq 2$.    Since $K\leq I_{G}(\Phi) \leq G$, $K/N$ contains the unique minimal normal subgroup of $G/N$ and $m(I_{G}(\Phi)/N)\leq 2$, we have that either   $I_{G}(\Phi)=G$ or $(G/N, I_{G}(\Phi)/N)\in \{(\SL(2,5),Z(\SL(2,5))), (C_2^4\rtimes A_7, C_2^4\rtimes A_6 ) \}$.

Let $T\subseteq \Irr(S^n)$ be the set of characters which are the product of $k$ copies of $\phi$ and $n-k$ copies of $1_S$. We have that $G$ acts on $T$ by conjugation. By the previous reasoning, we have that there are at most one orbit of size $1$ and for any of the two exceptions, there exists at most one orbit of size larger than $1$ (and its size is known). Since $|T|={n\choose k}={n\choose \lc n/2\rc}$, we have that ${n\choose \lc n/2\rc}=a+b \cdot c$ where $a,b \in \{0,1\}$ and $c\in \{7,60\}$. This forces $n=1$, which is impossible since $G/K\leq \Sigma_n$. 
\end{claim}

Thus, if $G$ is a group with $m(G)=2$, then $G$ possesses a unique non-abelian composition factor. Now, Theorem B follows from Theorems \ref{SolvableCase} and \ref{1cf}.

\renewcommand{\abstractname}{Acknowledgements}
\begin{abstract}
The author would like to thank Alexander Moret\'o for suggesting me this subject, and for many helpful conversations.  The author would also like to thank David Craven,  Derek Holt  and Gareth Tracey for their help with the computational issues of the work.
\end{abstract}



\end{document}